\documentclass[11pt,reqno,english]{amsart}
\usepackage[T1]{fontenc}
\usepackage[latin9]{inputenc}
\usepackage{geometry}
\usepackage{verbatim}
\usepackage{mathrsfs}
\usepackage{amsthm}
\usepackage{amstext}
\usepackage{amssymb}
\usepackage{graphicx}
\usepackage{esint}
\usepackage{color}
\usepackage{amsfonts}
\usepackage{amsmath} 
\usepackage{bm}           
\usepackage{mathabx}
\usepackage{amsxtra}

\usepackage{float,epsf}
\usepackage[english]{babel}
\usepackage{enumerate}
\usepackage{tikz}
\usepackage{srcltx}
\usepackage{verbatim}
\usepackage[bottom]{footmisc}
\makeatletter

\newcommand{\R}{\mathbb{R}}

\newcommand{\dd}{{\rm d\hspace{0.1mm}}}
\newcommand{\dr}{{\rm d}}

\newcommand{\M}{\mathcal{M}}

 \newcommand{\tn}{\textnormal}

 \newcommand{\Lip}{\text{\rm Lip}}
 
\theoremstyle{plain}

\newcommand{\FF}{{\boldsymbol F}}

\newcommand{\DM}{\mathcal D\mathcal M} 
\newcommand{\redb }{\partial^{*}} 
\newcommand{\eps}{\varepsilon}
\renewcommand{\div}{\text{\sl div}} 

\newcommand{\ii}{\rm i}

\newcommand{\ban}[1]{\left\langle  #1 \right\rangle}
\newcommand{\vv}{{\boldsymbol v}}

\newcommand{\bu}{{\bf u}}
\newcommand{\bbf}{{\bf f}}
\newcommand{\bq}{{\bf q}}

\newcommand{\res}{\mathop{\hbox{\vrule height 7pt width .5pt depth 0pt
\vrule height .5pt width 6pt depth 0pt}}\nolimits}

\newcommand{\Haus}[1]{{\mathscr H}^{#1}} 

\def\intave#1{\int_{#1}\hbox{\llap{$\raise2.3pt\hbox{\vrule
height.9pt width7pt}\phantom{\scriptstyle{#1}}\mkern-2mu$}}}

\def\intav#1{\mathchoice
          {\mathop{\vrule width 9pt height 3 pt depth -2.6pt
                  \kern -9pt \intop}\nolimits_{\kern -6pt#1}}%
          {\mathop{\vrule width 5pt height 3 pt depth -2.6pt
                  \kern -6pt \intop}\nolimits_{#1}}%
          {\mathop{\vrule width 5pt height 3 pt depth -2.6pt
                  \kern -6pt \intop}\nolimits_{#1}}%
          {\mathop{\vrule width 5pt height 3 pt depth -2.6pt
                  \kern -6pt \intop}\nolimits_{#1}}}
\def\intav#1{\vint_{#1}}

\allowdisplaybreaks[1]
\allowdisplaybreaks[1]

\usepackage{babel}

\usepackage{babel}

\makeatother

\usepackage{babel}

\begin{document}
\title[Divergence-Measure Fields: Gauss-Green Formulas and Normal Traces]{Divergence-Measure Fields:\\ Gauss-Green Formulas and Normal Traces}
\author{Gui-Qiang G. Chen}
\author{Monica Torres}
\maketitle

\bigskip
\smallskip
It is hard to imagine how most fields of science could provide the stunning mathematical descriptions
of their theories and results if the integration by parts formula did not exist.
Indeed,
integration by parts is an indispensable fundamental operation,
which has been used across scientific theories to pass from global (integral)
to local (differential) formulations of physical laws.
Even though the integration by parts formula
is commonly
known as the {\it Gauss-Green formula} (or the {\it divergence theorem}, or {\it Ostrogradsky's theorem}),
its discovery and rigorous mathematical proof are the result of the combined efforts of many
great
mathematicians,
starting back the period when the calculus was invented by Newton and Leibniz in the 17th century.

The one-dimensional integration by parts formula for smooth functions was first
discovered by
Taylor (1715).
The formula
is a consequence of the Leibniz product rule
and the Newton-Leibniz formula for the fundamental theorem of calculus.

The classical Gauss-Green formula for the multidimensional case is generally stated
for $C^{1}$ vector fields
and domains with $C^{1}$ boundaries.
However, motivated by the physical solutions with discontinuity/singularity
for Partial Differential Equations (PDEs) and Calculus of Variations,
such as nonlinear hyperbolic conservation laws and Euler-Lagrange equations,
the following fundamental issue arises:

\begin{quotation}
{\it Does the Gauss-Green formula still hold for vector fields with discontinuity/singularity
$($such as divergence-measure fields$)$ and domains with rough boundaries}?
\end{quotation}

The objective of this paper is to provide an answer to this issue
and to present a short historical review of the contributions by many mathematicians
spanning more than two centuries, which have made the discovery of the Gauss-Green formula possible.

\medskip
\section*{The Classical Gauss-Green Formula }

\begin{figure}[h]
\centering
\includegraphics[width=3.80in]{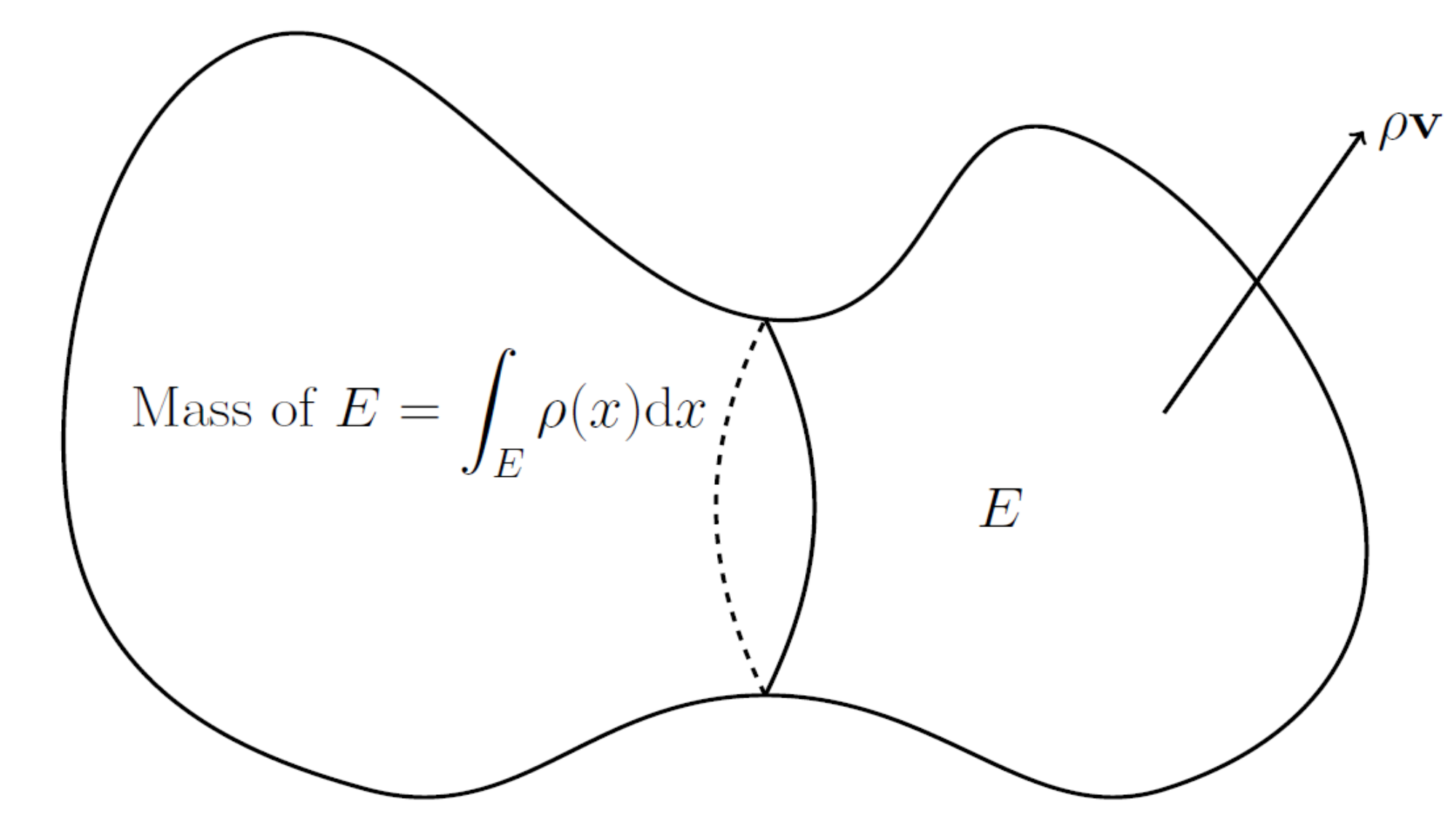}
\caption[]{Conservation of mass: The rate of change of the total mass in an open set $E$, $\frac{{\rm d}}{{\rm d}t}\int_E \rho(t,x)\,{\rm d}x$,
is equal to the total flux of mass across boundary $\partial E$, $\int_{\partial E}(\rho\mathbf{v})\cdot \nu\, \dr\Haus{n - 1}$,
where $\rho$ is the density, $\mathbf{v}$ is the velocity field, and the boundary value of $(\rho\mathbf{v})\cdot \nu$
is regarded as the normal trace of the vector field $\rho\mathbf{v}$ on $\partial E$.
The Gauss-Green formula yields the Euler equation for the conservation of mass:
$\rho_t+ \div(\rho\mathbf{v})=0$ in the smooth case.}\label{figure-1}
\end{figure}

The Gauss-Green formula was originally motivated in the analysis of fluids,
electric and magnetic fields, and other problems in the sciences
in order to establish the equivalence of integral and differential formulations of various physical
laws. In particular, the derivations of the Euler equations and the Navier-Stokes equations
in Fluid Dynamics and Gauss's laws
for the electronic and magnetic fields
are based on the validity of the Gauss-Green formula and associated Stokes theorem.
As an example, see Fig. \ref{figure-1} for the derivation of the Euler equation for the
conservation of mass in the smooth case.

\begin{figure}
\begin{minipage}{0.480\textwidth}
\centering
\includegraphics[height=1.4in,width=1.4in]{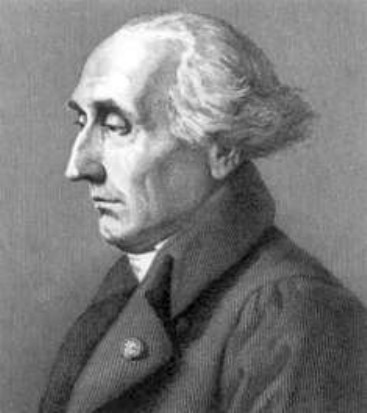}
\caption{\small Joseph-Louis Lagrange
(25 January 1736 -- 10 April 1813)}\label{figure-2}
\end{minipage}
\,\,
\begin{minipage}{0.480\textwidth}
\centering
\includegraphics[height=1.6in,width=2.3in]{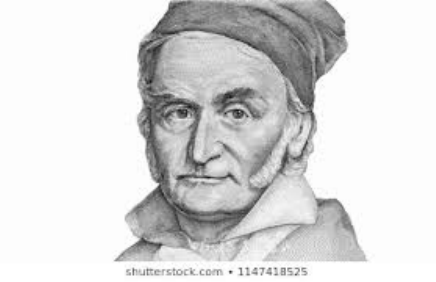}
\caption{\small Carl Friedrich Gauss (30 April 1777 -- 23 February 1855)}\label{figure-3}
\end{minipage}
\end{figure}

\begin{figure}
\centering
\centering
\includegraphics[height=1.61in,width=1.38in]{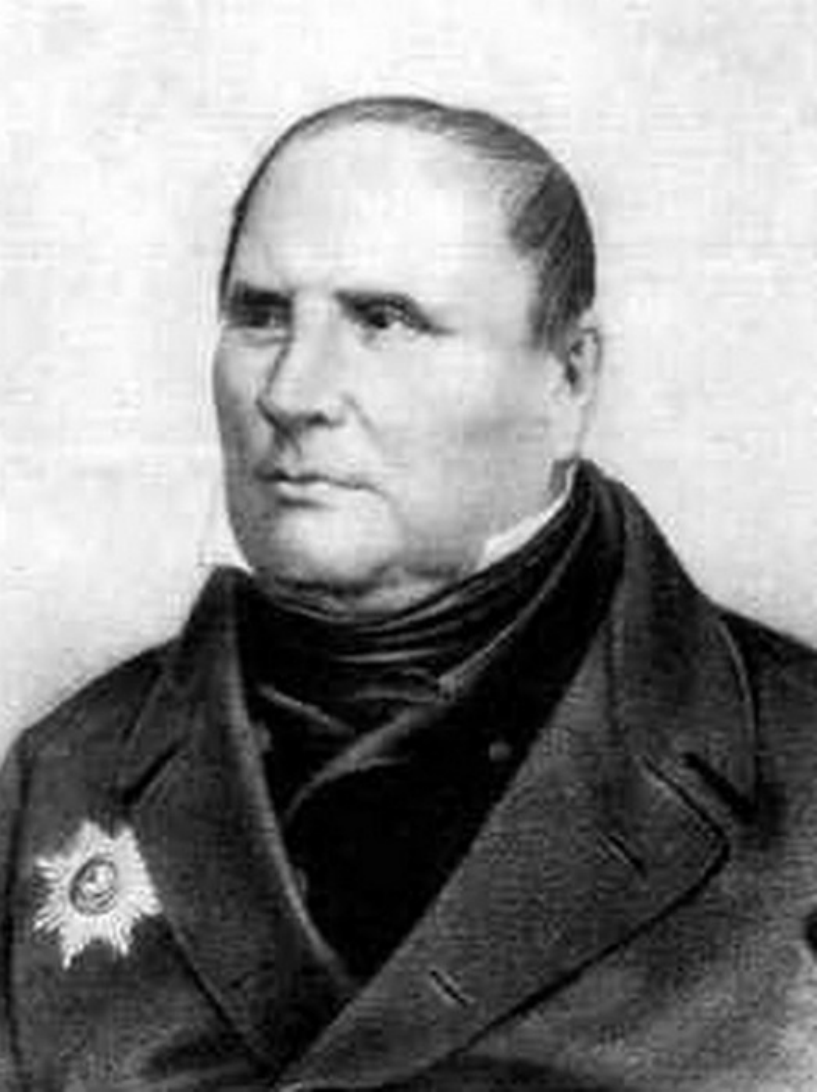}
\caption{\small Mikhail Ostrogradsky (24 September 1801 -- 1 January 1862)}\label{figure-4}
\end{figure}

\begin{figure}
\centering
\includegraphics[height=3.39in,width=2.91in]{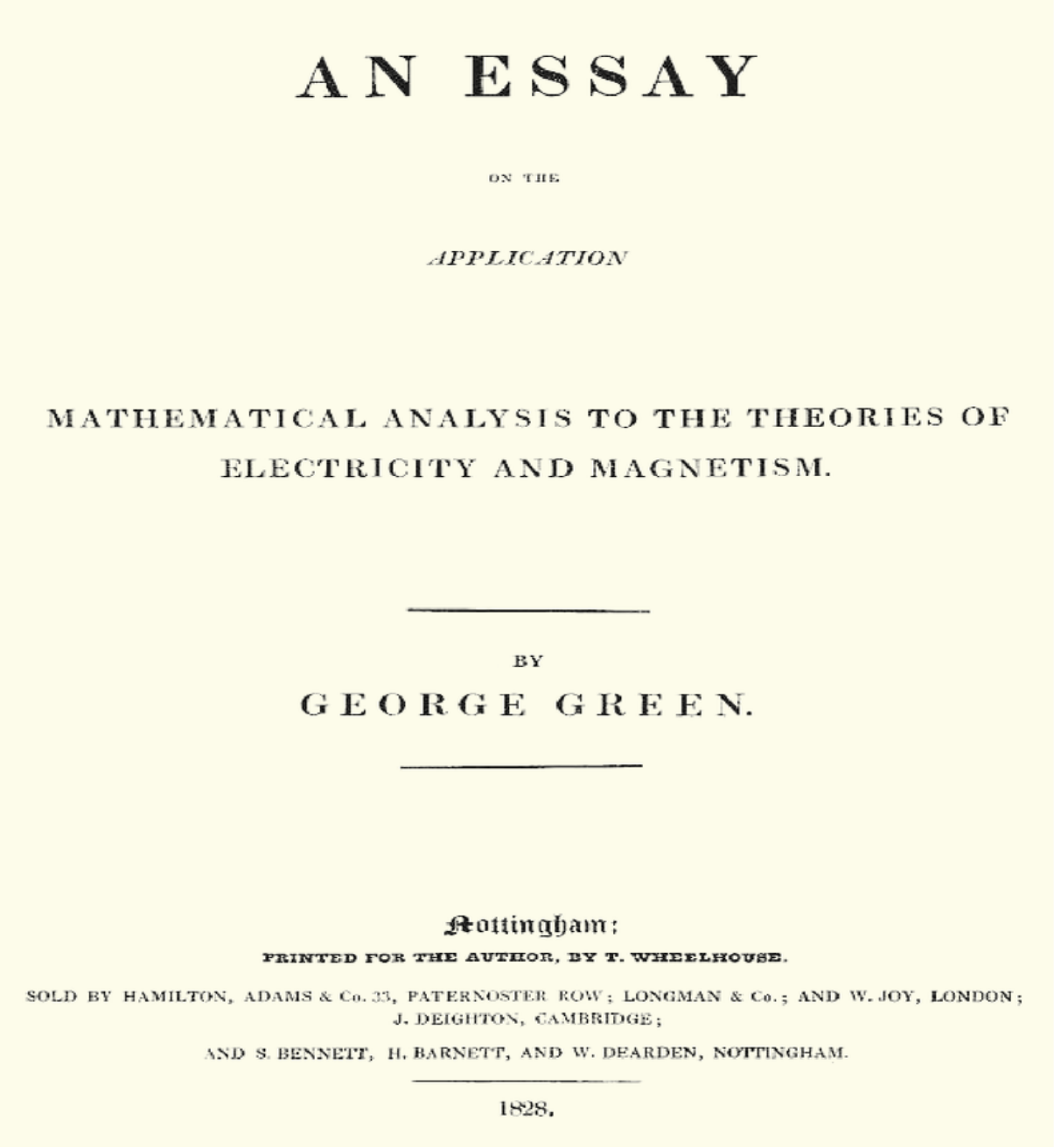}
\caption{\small The title page of the original essay of George Green (14 July 1793 -- 31 May 1841)}\label{figure-5}
\end{figure}

The formula that is also later known as the {\it divergence theorem} was first discovered by Lagrange\footnote{Lagrange, J.-L.:
{Nouvelles recherches sur la nature et la propagation du son},
{\it Miscellanea Taurinensia} (also known as: {\it M\'{e}langes de Turin}), 2: 11--172, 1762. He treated
a special case of the divergence theorem and transformed triple integrals into double integrals via integration by parts.}
in 1762 (see Fig. \ref{figure-2}),
however, he did not provide a proof of the result.
The theorem was later rediscovered by Gauss\footnote{Gauss, C.~F.: Theoria attractionis corporum sphaeroidicorum ellipticorum homogeneorum methodo nova tractata,
{\it Commentationes Societatis Regiae Scientiarium Gottingensis Recentiores}, 2: 355--378, 1813. In this paper,
a special case of the theorem was considered.} in 1813 (see Fig. \ref{figure-3})
and Ostrogradsky\footnote{Ostrogradsky, M. (presented on November 5, 1828; published in 1831):
Premi\`{e}re note sur la th\'{e}orie de la chaleur (First note on the theory of heat),
{\it M\'{e}moires de l'Acad\'{e}mie Imp\'{e}riale des Sciences de St. P\'{e}tersbourg},
Series 6, 1: 129--133, 1831. He stated and proved the {divergence-theorem} in
its Cartesian coordinate form.} in 1828 (see Fig. \ref{figure-4}).
Ortrogradsky's method of proof was similar to the approach Gauss used.
Independently, Green\footnote{Green, G.: {\it An Essay on the Application of Mathematical Analysis to the Theories of Electricity and Magnetism},
Nottingham, England: T. Wheelhouse, 1828.} (see Fig. \ref{figure-5})
also rediscovered the {\it divergence theorem} in the two-dimensional case
and published his result in 1828.

The {\it divergence theorem} in its vector form for the $n$--dimensional case $(n\ge 2)$
can be stated as
   \begin{equation}
  \label{formula}
  \int_{E} \div \FF \, \dr y= -\int_{\partial E}  \FF \cdot \nu  \, \dr\Haus{n - 1},
  \end{equation}
where $\FF$ is a $C^{1}$ vector field,  $E$ is a bounded open set with piecewise smooth boundary,
$\nu$ is the {\it inner} unit normal vector to $E$,
the boundary value of $\FF\cdot \nu$ is regarded as the {\it normal trace} of the vector field $\FF$
on $\partial E$,
and $\Haus{n-1}$ is the $(n-1)$--dimensional Hausdorff measure (that is an extension of the surface area
measure for $2$--dimensional surfaces
to general $(n-1)$--dimensional boundaries $\partial E$).
The formulation of \eqref{formula}, where $\FF$ represents a physical vector quantity,
is also the result of the efforts of many
mathematicians\footnote{See Stolze, C.~H.:  A history of the divergence theorem, In: {\it Historia Mathematica}, 437--442, 1978,
and the references therein.}
including
Gibbs, Heaviside, Poisson, Sarrus, Stokes, and Volterra.
In conclusion, formula \eqref{formula} is the result of more than two centuries of efforts by great mathematicians!

\medskip
\section*{Gauss-Green Formulas and Traces for Lipschitz Vector Fields \\ on Sets of Finite
Perimeter}

We first go back to the issue arisen earlier for extending the Gauss-Green formula to very rough sets.
The development of Geometric Measure Theory in the middle of the 20th century
opened the door to the extension
of the classical Gauss-Green formula over {\it sets of finite perimeter}
(whose boundaries can be very rough and contain cusps, corners, among others; {\it cf.} Fig. 6)
for Lipschitz vector fields.

Indeed, we may consider the left side of \eqref{formula} as a linear functional acting
on vector fields $\FF \in C_c^1(\mathbb{R}^n)$.
If $E$ is such that the functional: $\FF \to \int_{E} \div \FF\, \dr y$ is bounded
on $C_c(\R^n)$,
then the Riesz representation theorem implies that there exists a Radon measure $\mu_E$
such that
\begin{equation}\label{segunda}
\int_{E} \div \FF \, \dr y= \int_{\mathbb{R}^n} \FF \cdot \dr \mu_E \qquad\,\, \text{ for all } \FF \in C_c^1(\mathbb{R}^n),
\end{equation}
and the set, $E$, is called a {\it set of finite perimeter} in $\mathbb{R}^n$.
In this case, the Radon measure $\mu_E$ is actually the gradient of $-\chi_E$
(strictly speaking, the gradient is in the distributional sense),
where $\chi_E$ is
the characteristic function of $E$.
A set of density $\alpha\in [0,1]$ of $E$ in $\mathbb{R}^n$ is defined by
\begin{equation}\label{9a}
 E^{\alpha}: = \big\{y \in \R^{n}: \lim_{r \to 0} \frac {|B_r(y) \cap E| }{|B_r(y)|} = \alpha \big\},
 \end{equation}
where $|B|$ represents the Lebesgue measure of any Lebesgue measurable set $B$.
Then $E^0$ is the measure-theoretic exterior of $E$, while $E^1$ is the measure-theoretic interior of $E$.

\vspace{-10pt}
\begin{figure}[h]
\centering
\includegraphics[width=4.2in]{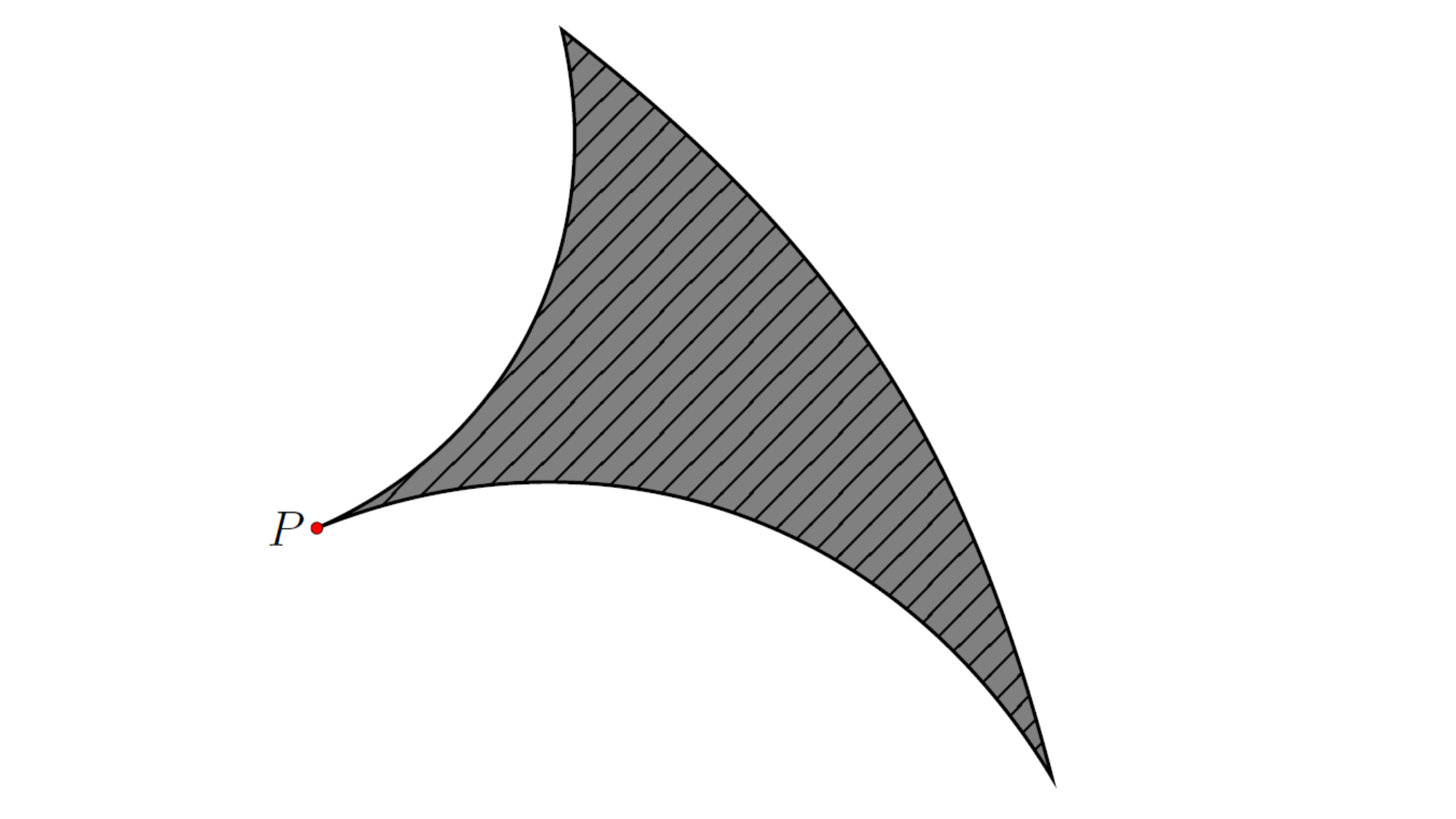}

\vspace{-15pt}
\caption[]{$E$ is an open set of finite perimeter where cusp $P$ is a point of density $0$}
\end{figure}

Sets of finite perimeter can be quite subtle.
 For example, the countable union
 of open balls with centers on the rational points $y_k$, $k=1,2,\dots$, of the unit ball in $\R^{n}$ and
 with radius $2^{-k}$, is a set of finite perimeter.
 A set of finite perimeter may have
 a large set of {\it cusps} in the topological boundary ({\it e.g.}, $\Haus{n-1}(E^{0} \cap \partial E)> 0$ or $\Haus{n-1}(E^{1} \cap \partial E)> 0$).
The set $E$ can also have points in the boundary which belong $E^{\alpha}$
 for $\alpha \notin \{0,1\}$.
For example, the four corners of a square are points of density $\frac{1}{4}$.

Even though the boundary of a set $E$ of finite perimeter can be very rough,
De Giorgi's structure theorem
indicates that it has nice tangential properties so that
there is a notion of measure-theoretic tangent plane.
More rigorously, the topological boundary $\partial E$ of $E$
contains an $(n-1)$--rectifiable set,
known as the reduced boundary of $E$,  denoted as $\partial^* E$,
which can be covered by a countable union of $C^1$ surfaces, up to a set of $\Haus{n-1}$--measure zero.
More generally, Federer's theorem states that $\Haus{n-1} (\partial^s E \setminus \partial^* E) =0$,
 where $\partial^{s} E=\R^{n} \setminus (E^{0} \cup E^{1})$ is called the {\it essential boundary}.
Note that $\partial^* E \subset E^{\frac{1}{2}}$ and $\partial^* E \subset \partial^s E \subset \partial E$.

It can be shown that every $y \in \partial^* E$ has an {\it inner} unit normal vector $\nu_{E}(y)$
and a tangent plane in the measure-theoretic sense,
and
\eqref{segunda} reduces to
\begin{equation}
\label{segunda2}
\int_E \div \FF \, \dr y = -\int_{\partial^* E} \FF (y) \cdot \nu_E (y) \, \dr \mathscr{H}^{n-1}(y).
\end{equation}
This Gauss-Green formula for Lipschitz vector fields $\FF$ over sets of finite perimeter
was proved by De Giorgi (1954--55)
and Federer (1945, 1958)
in a series of papers.
See \cite{Federer,Pfeffer} and the references therein.

\section*{Gauss-Green Formulas and Traces for Sobolev and $BV$ Functions \\
on Lipschitz Domains}
It happens in many areas of analysis, such as PDEs and Calculus of Variations,
that it is necessary to work with the functions that are not Lipschitz, but only in $L^{p}$, $1 \leq p \leq \infty$.
In many of these cases, the functions have
weak
derivatives that also belong to $L^{p}$.
That is, the corresponding $\FF$ in \eqref{segunda2} is a Sobolev vector field.
The necessary and sufficient conditions for the existence of traces ({\it i.e.}, boundary values)
of Sobolev functions defined on the boundary of the domain
have been obtained so that \eqref{segunda2} is a valid formula over open sets with Lipschitz boundary.

The development of the theory of Sobolev spaces has been fundamental in analysis.
However, for many further applications, this theory is still not sufficient.
For example, the characteristic function of a set $E$ of finite perimeter, $\chi_E$, is not a Sobolev function.
Physical solutions in gas dynamics involve shock waves and vortex sheets that are discontinuities with jumps.
Thus, a larger space of functions, called the space of functions of bounded variation ($BV$), is necessary,
which consists of all functions in $L^1$ whose
derivatives are Radon measures\footnote{The definition of the $BV$ space is
a generalization of the classical notion of the class of one-dimensional functions $f(x)$ with finite total variation (TV)
over an interval $[a,b]\subset \R${\rm :}
$
TV_a^b(f)=\sum_{\mathcal{P}}|f(x_{j+1})-f(x_j)|,
$
where the supremum runs over the set of all partitions
$\mathcal{P}=\{(x_0,\cdots,x_J)\,: \, a\le x_0\le\cdots\le x_J\le b, 1\le J<\infty\}$.
}.

This space has compactness properties that allow, for instance, to show the existence of minimal surfaces
and the well-posedness of $BV$ solutions for hyperbolic conservation laws.
Moreover, the Gauss-Green formula \eqref{segunda2} is also valid for $BV$ vector fields over Lipschitz domains.
See \cite{Federer,maz2013sobolev,VH} and the references therein.

\section*{Divergence-Measure Fields and Hyperbolic Conservation Laws}

A vector field $\FF\in L^p (\Omega)$, $1\le p\le \infty$, is called
a divergence-measure field
if $\div \FF$ is a signed Radon  measure
with finite total variation in $\Omega$.
Such vector fields form Banach spaces, denoted as $\mathcal{DM}^p(\Omega)$, for  $1\le p\le \infty$.

These spaces arise naturally in the field of Hyperbolic Conservation Laws.
Consider hyperbolic systems of conservation laws of the form:
\begin{eqnarray}
\label{codo}
\bu_t+ \nabla_x \cdot \bbf(\bu)=0 \qquad\,\,\,\,\mbox{for $\bu=(u_1,...,u_m)^\top: \R^{n} \to \R^m$}
\end{eqnarray}
where $(t,x)\in \R_+^n:=[0,\infty)\times \R^d$, $n:=d+1, \bbf=(\bbf_1,\bbf_2,...,\bbf_d)$, and $\bbf_j:\R^m \to \R^m, j=1,\cdots,d$.
A prototype of such systems is the system of Euler equations for compressible fluids,
which consists of the conservation equations of mass, momentum, and energy in Continuum Mechanics
such as Gas Dynamics and Elasticity.
One of the main features of the hyperbolic system is that the speeds of propagation of solutions are finite;
another feature is that, no matter how smooth a solution starts with initially, it generically develops
singularity and becomes a discontinuous/singular solution.  To single out physically relevant solutions,
it requires the solutions of system \eqref{codo} to fall within the following class
of entropy solutions:
\begin{quotation}
{\it An entropy solution $\bu(t,x)$
of system \eqref{codo} is characterized by the Lax entropy inequality{\rm :}
For any entropy pair $(\eta, \bq)$ with $\eta$ being a convex function of $\bu$,
\begin{equation}\label{6-a}
\eta(\bu(t,x))_t + \nabla_x\cdot \bq(\bu(t,x)) \leq 0 \qquad \tn{\it in the sense of distributions}.
\end{equation}
}
\end{quotation}
Here, a function $\eta\in C^1(\R^m,\R)$
is called an entropy of system \eqref{codo}
if there exists an entropy flux
$\bq=(q_1,\cdots, q_d)\in C^1(\R^m, \R^d)$ such that
\begin{equation}
\label{5-a}
\nabla q_j(\bu)= \nabla \eta (\bu)  \nabla \bbf_j (\bu)  \qquad\,\, \mbox{for $j=1,2,...,d$}.
\end{equation}
Then $(\eta, \bq)(\bu)$ is called an entropy pair of system \eqref{codo}.

Friedrichs-Lax (1971) observed that most systems of conservation
laws that result from Continuum Mechanics are endowed with a globally
defined, strictly convex entropy.
In particular, for the compressible Euler equations in Lagrangian coordinates,
$\eta=-S$ is such an entropy, where $S$ is the physical thermodynamic entropy of the fluid.
Then the Lax entropy inequality for the physical entropy $\eta=-S$ is an exact statement of
the Second Law of Thermodynamics
({\it cf}. \cite{Chen2,dafermos2010hyperbolic}).
Similar notions of entropy have also been used in many
fields such as Kinetic Theory, Statistical Physics, Ergodic Theory, Information Theory,
and Stochastic Analysis.

Indeed, the available existence theories show that the solutions of \eqref{codo}
are entropy solutions obeying the Lax entropy inequality \eqref{6-a}.
This implies that, for any entropy pair $(\eta, \bq)$ with $\eta$ being a convex function,
there exists  a nonnegative measure $\mu_{\eta} \in \M(\R_+^n)$ such that
\begin{equation}
\label{7-a}
\div_{(t,x)} (\eta(\bu(t,x)), \bq(\bu(t,x)))= -\mu_{\eta}.
\end{equation}
Moreover, for any $L^\infty$ entropy solution $\bu$,
if the system is endowed with a strictly
convex entropy, then, for any $C^2$ entropy pair $(\eta, \bq)$ (not necessarily convex for $\eta$),
there
exists a signed Radon measure $\mu_\eta\in \mathcal{M}(\R_+^n)$ such that \eqref{7-a}
still holds (see \cite{Chen2}).
For these cases,  $(\eta(\bu), \bq(\bu))(t,x)$ is a $\mathcal{DM}^p(\R_+^n)$ vector field,
as long as $(\eta(\bu), \bq(\bu))(t,x)\in L^p(\R_+^n, \R^n)$ for some $p\in [1,\infty]$.

Equation \eqref{7-a} is one of the main motivations to develop a $\DM$ theory in  Chen-Frid \cite{CF1,CF2}.
In particular, one of the major issues is whether integration by parts can be performed  in \eqref{6-a}
to explore to fullest extent possible all the information about the entropy solution $\bu$.
Thus, a concept of normal traces for $\DM$ fields $\FF$ is necessary to be developed.
The existence of normal traces is also fundamental for initial-boundary value problems for
hyperbolic systems \eqref{codo}
and for the analysis of structure and regularity of entropy solutions $\bu$
(see {\it e.g.}  \cite{Chen2,CF1,CF2,dafermos2010hyperbolic,Vas}).

Motivated by hyperbolic conservation laws,
the {\it interior} and {\it exterior} normal traces need to be constructed as the limit of classical normal traces on one-sided smooth approximations of the domain.
Then the surface of a {\it shock wave} or {\it vortex sheet} can be approximated with {\it smooth surfaces}
to obtain the {\it interior} and {\it exterior} fluxes on the {\it shock wave} or {\it vortex sheet}.

Other important connections for the $\DM$ theory
are the characterization of phase transitions (coexistent phases with discontinuities across the boundaries),
the concept of stress, the notion of Cauchy flux, and the principle of balance law
to accommodate the discontinuities and singularities in the continuum media
in Continuum Mechanics, as discussed
in Degiovanni-Marzocchi-Musesti \cite{DMM}, \v{S}ilhav\'y \cite{Silhavy1,Silhavy1b},
Chen-Torres-Ziemer \cite{ctz}, and Chen-Comi-Torres \cite{cct}.

\section*{Gauss-Green Formulas and Normal Traces for $\mathcal{DM}^{\infty}$ Fields}

We start with the following simple example:

\medskip
\noindent{\bf Example 1}. {\it
Consider the vector field $\FF: \Omega=\R^{2} \cap  \{y_1>0\} \to \R^2${\rm :}
 \begin{equation}
 \FF(y_1,y_2) =(0,  \sin(\frac{1}{y_1})).
 \end{equation}
Then $\FF \in \DM^\infty(\Omega)$ with $\div \FF=0$ in $\Omega$.

However, in an open half-disk $E=\{(y_1,y_2)\,:\, y_1>0, y_1^2+y_2^2<1\}
\subset \Omega$,
the previous Gauss-Green formulas do not apply.
Indeed, since $\FF(y_1,y_2)$ is highly oscillatory when $y_1\to 0$ which is not well-defined at $y_1=0$,
it is not clear how the normal trace
$\FF \cdot \nu$
on $\{y_1=0\}$ can be understood in the classical sense
so that the equality between $\int_{E} \div \FF\,\dr\Haus{n-1} =0$
and $\int_{\partial E} \FF \cdot \nu\,\dr\Haus{n-1}$ holds.
This example shows that a suitable notion of normal traces is required to be developed.
}

\medskip
A generalization of \eqref{formula} to $\DM^\infty(\Omega)$  fields
and bounded sets with Lipschitz boundary was derived in
Anzellotti \cite{Anzellotti_1983} and Chen-Frid \cite{CF1} by different approaches.
A further generalization of \eqref{formula} to $\DM^\infty(\Omega)$ fields
and arbitrary bounded sets of finite perimeter, $E \Subset \Omega$, was first obtained
in Chen-Torres \cite{ChenTorres} and \v{S}ilhav\'y \cite{Silhavy1} independently;
see also \cite{ctz,ComiTorres}.

\medskip
\noindent{\bf Theorem 1}.
{\it Let $E$ be a set of finite perimeter. Then
 \begin{equation}
 \label{gen1}
  \int_{E^{1}} \phi\, \dd \div \FF + \int_{E^1} \FF \cdot \nabla \phi \,\dr y
   = - \int_{\redb E} \phi\, \mathfrak{F}_{\ii} \cdot \nu_{E} \, \dr \Haus{n - 1}
 \qquad\mbox{for every $\phi \in C_c^{1}(\R^n)$},
 \end{equation}
 where the interior
 normal trace $\mathfrak{F}_{\ii} \cdot \nu_{E}$
 is a bounded function defined on the reduced boundary of $E$
 {\rm (}{\it i.e.}, $\mathfrak{F}_{\ii} \cdot \nu_{E}\in L^{\infty}(\redb E; \Haus{n - 1})${\rm )},
 and $E^{1}$ is the measure-theoretic interior of $E$ as defined in \eqref{9a}.}

\smallskip
 One approach for the proof of \eqref{gen1} is based on a product rule for $\DM^\infty$
 fields (see \cite{ChenTorres}).
 Another approach in \cite{ctz}, following \cite{CF1}, is based on a new approximation theorem for sets of finite perimeter,
 which shows that the level sets of convolutions $w_k:=\chi_E * \varrho_k$ by the standard positive and symmetric mollifiers
 provide smooth approximations essentially
 from the {\it interior} (by choosing  $w_k^{-1}(t)$ for $\frac{1}{2}< t <1$)
 and the {\it exterior} (for $0 <t < \frac{1}{2}$).
 Thus, the interior normal trace $\mathfrak{F}_{\ii} \cdot \nu_{E}$
 is constructed as the limit of classical normal traces over the smooth approximations
 of $E$.
 Since the level set
 $w_k^{-1}(t)$ (with a suitable fixed $0 < t< \frac{1}{2}$)
 can intersect the measure-theoretic exterior $E^0$ of $E$,
 a critical step for this approach
 is to show that $\Haus{n-1} (w_k^{-1}(t) \cap E^{0})$ converges to zero
 as $k \to \infty$.
 A key point for 
 this proof is the fact that, if $\FF$ is a $\DM^\infty$ field,
 then the Radon measure $|\div \FF|$ is absolutely continuous with respect to $\Haus{n-1}$, as first observed by Chen-Frid \cite{CF1}.

\smallskip
If $E=\{y \in \R^{2}\,:\, |y|<1\}\setminus\{y_1>0, y_2=0\}$,
the above formulas apply, but
the integration is not over the original representative
consisting of the disk with radius $\{y_1>0, y_2=0\}$ removed, since $E^{1}= \{y \in \R^{2}\,:\,|y|<1 \}$
is the open disk.
In many applications including those in Materials Science, we may want to integrate on a domain with {\it fractures} or {\it cracks}; see also Fig. 7.
Since the cracks are part of the topological boundary and belong to the measure-theoretic interior $E^{1}$,
the formulas in \eqref{gen1}
do not provide such information.
In order to establish a Gauss-Green formula that includes such cases,
we restrict to {\it open sets $E$ of finite perimeter} with $\Haus{n-1}(\partial E \setminus E^0) < \infty$.
Therefore, $\partial E$ can still have a large set of {\it cusps} or points
of density $0$ ({\it i.e.}, points belonging to $E^0$).

\begin{figure}[h]
\centering
\includegraphics[width=3.8in]{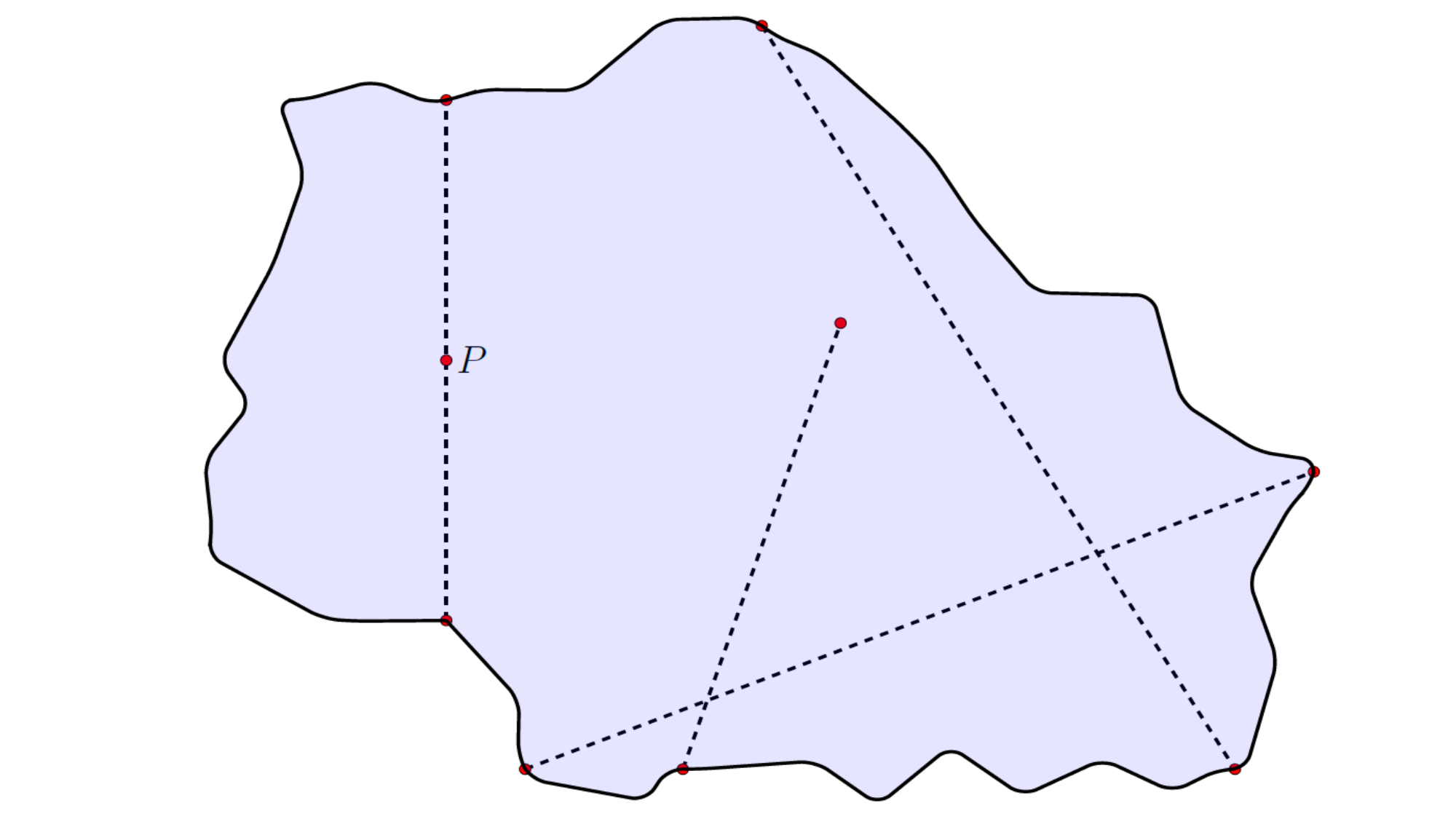}

\vspace{-10pt}
\caption[]{$E$ is an open set with several fractures (dotted lines)}
\end{figure}

It has been shown in \cite{CLT} that, if $E$ is a bounded open set satisfying
\begin{equation}\label{11a}
\Haus{n-1}(\partial E \setminus E^0) < \infty,
\end{equation}
and $\FF\in \mathcal{DM}^\infty(E)$,
then there exists a family of  sets $E_k \Subset E$ such that $E_k \to E$ in $L^1$ and $\sup_{k} \Haus{n-1}(\partial^* E_k) < \infty$
such that,
for any $\phi \in C_c^{1}(\R^n)$,
\begin{equation}
\label{nuevo}
\int_{E}\phi  \, \dd \div \FF
+ \int_{E} \FF \cdot \nabla \phi \,\dr y
=- \lim_{k \to \infty} \int_{\partial E_{k}} \phi \FF \cdot \nu_{E_{k}} \,\, \dr\Haus{n - 1}
=
\int_{\partial E \setminus E^0} \phi\,\mathfrak{F}\cdot\nu_E\, \dr \Haus{n-1},
\end{equation}
where $\mathfrak{F}\cdot\nu_E$ is well-defined in $L^{\infty}(\partial E\setminus E^0; \Haus{n - 1})$ as the interior normal trace on $\partial E\setminus E^0$.

This approximation result for the bounded open set $E$
can be accomplished by performing
delicate covering arguments, especially by applying the Besicovitch theorem
to a covering of $\partial E \cap E^{0}$.
Moreover, \eqref{nuevo} is a formula {\it up to the boundary},
since we do not assume that the domain of integration is compactly contained in the domain of $\FF$.
More general product rules for $\div (\phi \FF)$ can be proved to weaken
the regularity of $\phi$; see \cite{cct,CLT} and the references therein.

\section*{Gauss-Green Formulas and Normal Traces for $\mathcal{DM}^{p}$ Fields}
For $\DM^p$  fields with $1\le p<\infty$,
the situation becomes more delicate.

\begin{figure}
\centering
\begin{minipage}{0.495\textwidth}
\centering
\includegraphics[height=1.95in,width=3.4in]{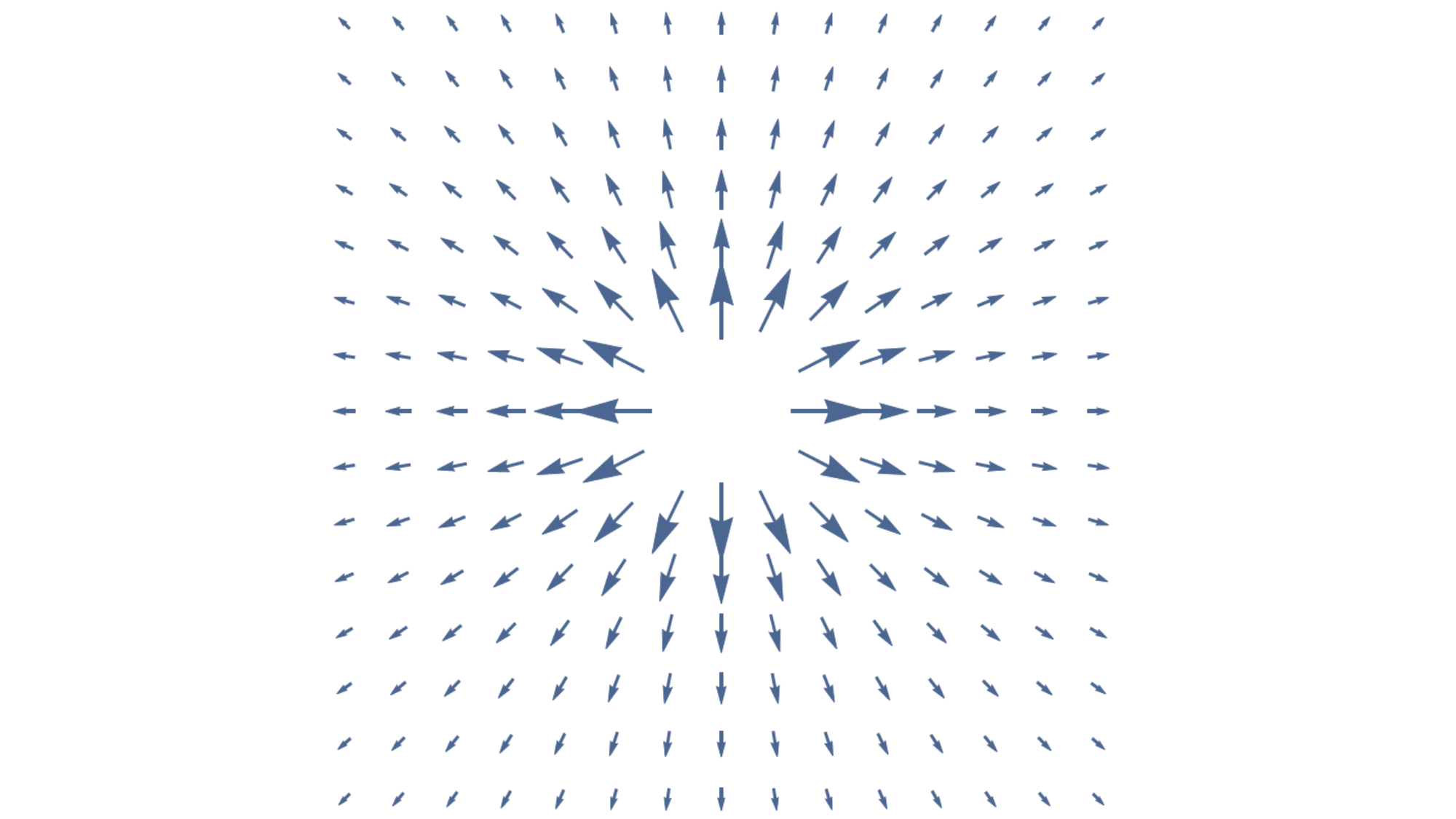}
\caption{\small The vector field $\FF(y_1,y_2)=\frac{(y_1,y_2)}{y_1^2+y_2^2}$ in Example 2}\label{figure-8}
\end{minipage}
\begin{minipage}{0.495\textwidth}
\centering
\includegraphics[height=1.7in,width=3.1in]{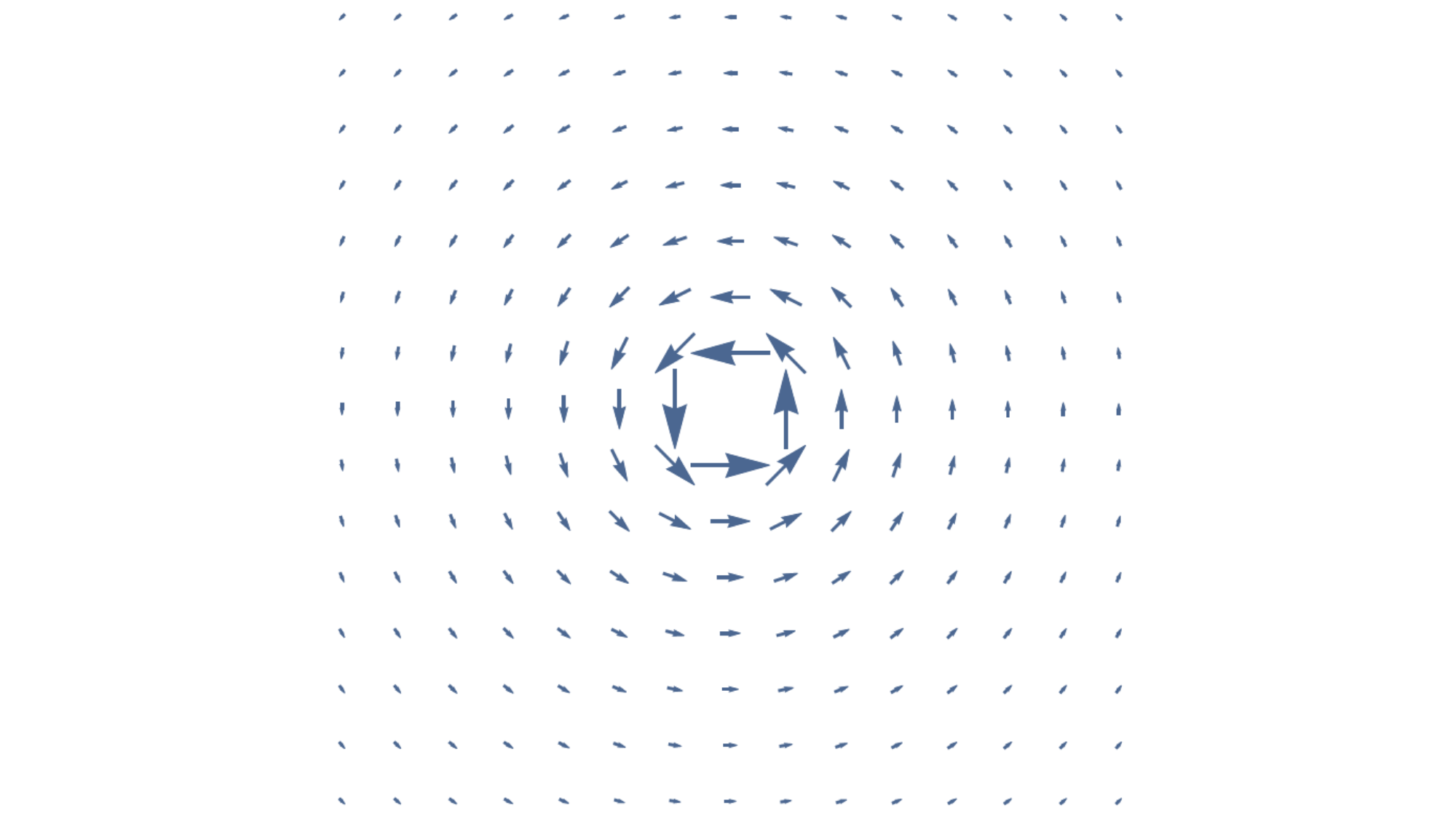}
\caption{\small The vector field $\FF(y_1,y_2)=\frac{(-y_2, y_1)}{y_1^2+y_2^2}$ in Example 3}\label{figure-9}
\end{minipage}
\end{figure}

\medskip
\noindent
{\bf Example 2} (Whitney [Example 1, p. 100]\footnote{Whitney, H.:
Geometric Integration Theory, Princeton University Press, Princeton, 1957.}).
{\it Consider the vector field  $\FF : \R^{2} \setminus \{ (0, 0) \} \to \R^{2}$ $($see Fig. {\rm \ref{figure-8}}$)${\rm :}
\begin{equation} \label{Whitney_field}
\FF(y_1, y_2)= \frac{(y_1,y_2)}{y_1^{2} + y_2^{2}}.
\end{equation}
Then $\FF \in \DM^{p}_{\rm loc}(\R^{2})$ for $1 \le p < 2$.
If $E = (0, 1)^{2}$, it is observed
that
\begin{equation*}
0 = \div \FF (E) \neq - \int_{\partial E} \FF \cdot \nu_{E} \, \dr \Haus{1} = \frac{\pi}{2},
\end{equation*}
where $\nu_{E}$ is the inner unit normal to the square.
However, if
$E_{\eps} := \{ y \in E\,:\,{\rm dist}(y, \partial E)> \eps \}$ and  $E^{\eps} := \{ y \in \R^{n}\,:\, {\rm dist}(y, E)<\eps \}$  for any $\eps > 0$,
then
\begin{align*}
&0 = \div \FF (E) = - \lim_{\eps \to 0} \int_{\partial E_{\eps}} \FF \cdot \nu_{E_{\eps}} \, \dr \Haus{1}, \\
&2 \pi = \div \FF (\overline{E})  = - \lim_{\eps \to 0} \int_{\partial E^{\eps}} \FF \cdot \nu_{E^{\eps}} \, \dr \Haus{1}.
\end{align*}
In this sense, the equality is achieved on both sides of the formula.
}

Indeed,
for a $\DM^p$ field with $p \neq \infty$,
$|\div \FF|$ is absolutely continuous with respect to $\Haus{n-p'}$ for $p'>1$
with $\frac{1}{p'}+\frac{1}{p}=1$, but not with respect to $\Haus{n-1}$ in general.
This implies that the approach in \cite{ctz}
does not apply directly to obtain normal traces for $\mathcal{DM}^p$ fields
for $p \neq \infty$.

Then the following questions arise:

\begin{itemize}

\item {\it Can the previous formulas be proved in general for any $\FF \in \mathcal{DM}^p(\Omega)$ and for any open set $E \subset \Omega$}?

\item {\it Even though almost all the level sets of the distance function are sets of finite perimeter,
can the formulas with smooth approximations of $E$ be obtained, in place of $E_{\eps}$ and $E^{\eps}$}?

\item {\it If $E$ has a Lipschitz boundary, do {\it regular Lipschitz deformations} of $E$, as defined in Chen-Frid \cite{CF1,CF2}, exist}?
\end{itemize}

The answer to all three questions is {\it affirmative}.

\medskip
\noindent
{\bf Theorem 2} (Chen-Comi-Torres \cite{cct}).
{\it Let $E \subset \Omega$ be a bounded open set, and let $\FF \in \DM^p(\Omega)$ for $1\le p\le \infty$.
Then, for any $\phi \in (C^{0}\cap L^{\infty})(\Omega)$ with $\nabla \phi \in L^{p'}(\Omega; \R^{n})$ for $p'\ge 1$ so that $\frac{1}{p'}+\frac{1}{p}=1$,
there exists a sequence of bounded open sets $E_{k}$ with $C^\infty$ boundary such that $E_{k} \Subset E$, $\bigcup_{k} E_{k} = E$, and
\begin{equation}
\label{smooth1}
\int_{E} \phi \, \dd \div \FF + \int_{E} \FF \cdot \nabla \phi \, \dr y
= - \lim_{k \to \infty} \int_{\partial E_{k}} \phi \FF \cdot \nu_{E_{k}} \,\, \dr\Haus{n - 1},
\end{equation}
where $\nu_{E_{k}}$ is the classical inner unit normal vector to $E_{k}$.}

\smallskip
For the open set $E$ with Lipschitz boundary,
it can be proved
that the deformations of $E$ obtained with the method of regularized distance are bi-Lipschitz.
We can also employ an alternative construction  by Ne{\v{c}}as (1962)
to obtain smooth approximations $E_{\eps}$ of a bounded open set $E$  with Lipschitz boundary in such a way
that the deformation, $\Psi_{\eps}(x)$, mapping $\partial E$ to $\partial E_{\eps}$ is bi-Lipschitz,
and the Jacobians of the deformations, $J^{\partial E}(\Psi_{\eps})$, converge to $1$ in $L^{1}(\partial E)$ as $\eps$ approaches
zero (see \cite[Theorem 8.19]{cct})\footnote{A simpler construction of smooth approximations $E_{\eps}$ can also be
obtained when $E$ is a {\it strongly} Lipschitz domain -- a stronger requirement than a general Lipschitz domain;
see Hofmann, S., Mitrea, M., Taylor, M.:
Geometric and transformational properties of Lipschitz domains, Semmes-Kenig-Toro domains, and other classes
of finite perimeter domains, {\it J. Geometric Anal.} 17: 593--647, 2007.}.
This shows that any bounded open set with Lipschitz boundary
admits a  {\it regular Lipschitz deformation}  in the sense of Chen-Frid \cite{CF1,CF2}.
Therefore, we can write more explicit Gauss-Green formulas for Lipschitz domains.

\medskip
\noindent
{\bf Theorem 3} (Chen-Comi-Torres \cite{cct}).
{\it If $E \Subset \Omega$ is an open set with Lipschitz boundary and $\FF \in \DM^p(\Omega)$ for $1 \le p \le \infty$,
then, for every $\phi \in C^{0}(\Omega)$ with $\nabla \phi \in L^{p'}(\Omega; \R^{n})$,
there exists a set $\mathcal{N} \subset \R$ of Lebesgue measure zero
such that,
for every nonnegative sequence $\{\eps_k\}\not\subset\mathcal{N}$ satisfying $\eps_k \to 0$,
\begin{equation*} \label{G-G area eps int}
\int_{E} \phi \, \dd \div \FF + \int_{E} \FF \cdot \nabla \phi \, \dr y
= - \lim_{k \to \infty} \int_{\partial E}
\big ( \phi \FF \cdot \frac{\nabla \rho}{|\nabla \rho|} \big) (\Psi_{\eps_{k}}(y)) J^{\partial E}(\Psi_{\eps_k})(y) \, \dr \Haus{n - 1}(y),
\end{equation*}
where $\rho(y)\in C^2(\R^n\setminus\partial E)\cap {\rm Lip(\R^n)}$ is a regularized
distance\footnote{Such a regularized distance has been constructed in Ball-Zarnescu \cite[Proposition 3.1]{Ball-Zarnescu};
also see Lieberman, G. M.: Regularized distance and its applications,
{\it Pac. J. Math.} 117(2): 329--352, 1985, and Fraenkel, L.~E.: On regularized distance and related functions, 
{\it Proc. R. Soc. Edinb. Sect. A Math.} 83: 115--122 (1979).} for $E$,
so that the ratio functions $\frac{\rho(y)}{d(y)}$ and $\frac{d(y)}{\rho(y)}$ are positive 
and uniformly bounded for all $y\in \R^n\setminus\partial U$
for $d(y)$ to be ${\rm dist}(y,\partial E)$ for $y\in E$ and
$-{\rm dist}(y,\partial E)$ for $y\notin E$.
}

\smallskip
The question now arises as to
whether
the limit can be realized on the right hand side of the previous formulas as an integral on $\partial E$.
In general, this is not possible.
However, in some cases, it is possible to represent the normal trace with a measure supported on $\partial E$.
In order to see this,  for $\FF \in \DM^{p}(\Omega)$,  $1 \le p \le \infty$,
and a bounded Borel set $E \subset \Omega$,  we
follow \cite{Anzellotti_1983,CF2,Silhavy1} to define the normal trace distribution of $\FF$ on $\partial E$ as
\begin{equation} \label{normal trace def}
\ban{\FF \cdot \nu, \phi}_{\partial E} := \int_{E} \phi \, \dr \div \FF + \int_{E} \FF \cdot \nabla \phi \, \dr y
\qquad \mbox{for any $\phi \in \Lip_{c}(\R^{n})$}.
\end{equation}
The formula presented above shows that the trace distribution on $\partial E$ can be extended to a functional
on $\{ \phi \in (C^{0} \cap L^{\infty})(\Omega): \nabla \phi \in L^{p'}(\Omega, \R^{n})\}$ so that
we can always represent the normal trace distribution as the limit of classical normal traces
on smooth approximations of $E$.
Then the question is whether there exists a Radon measure $\mu$ concentrated on $\partial E$ such that
$\ban{\FF\cdot \nu, \phi}_{\partial E}= \int_{\partial E} \phi\, \dr \mu$.
Unfortunately, this not the case in general.

\medskip
\noindent
{\bf Example 3}. Consider the vector field
$\FF: \R^2\setminus\{(0,0)\}\to \R^2$:
$$
\FF(y_1,y_2)=\frac{(-y_2, y_1)}{(y_1^2+y_2^2)^{\alpha/2}} \qquad \mbox{for $1\le \alpha <3$}.
$$
Then $\FF\in\mathcal{DM}_{\rm loc}^p(\R^2)$ with $1\le p<\frac{2}{\alpha-1}$ for $1<\alpha<3$, and
$p=\infty$ for $\alpha=1$, and $\div \FF=0$ on $E=(-1,1)\times (-1,0)$.
For $\phi\in {\rm Lip}(\partial E)$ with ${\rm supp}\,\phi\subset \{y_2=0, |y_1|<1\}$,
as shown in \v{S}ilhav\'y \cite[Example 2.5]{Silhavy2},
$$
\ban{\FF \cdot \nu, \phi}_{\partial E} :=
\begin{cases}
\int_{-1}^1 \phi(t,0){\rm sgn}(t)|t|^{1-\alpha} \dd t  \qquad&\mbox{for $1\le \alpha<2$},\\[1.5mm]
\lim_{\eps \to 0}\int_{\{|t|>\eps\}}
\phi(t,0){\rm sgn}(t)|t|^{1-\alpha} \dd t  \qquad&\mbox{for $2\le \alpha<3$}.
\end{cases}
$$
This indicates that the normal trace $\ban{\FF\cdot \nu, \cdot}_{\partial E}$
of the vector field $\FF(y_1,y_2)$
is a measure when $1\le \alpha<2$, but is not a measure when $2\le \alpha<3$.
Also see  Fig. \ref{figure-9} when $\alpha=2$.

\medskip
Indeed, it can be shown (see \cite[Theorem 4.1]{cct})
that  $\ban{\FF\cdot \nu, \phi}_{\partial E}$
can be represented as a measure if and only if   $\chi_{E} \FF \in \DM^{p}(\Omega)$.
Moreover, if  $\ban{\FF \cdot \nu, \cdot}_{\partial E}$ is a measure,
then
\begin{itemize}
\item[\rm (i)] For $p = \infty$,
$|\ban{\FF\cdot \nu, \cdot}_{\partial E}| \ll \Haus{n - 1} \res \partial E$
({\it i.e.}, $\Haus{n - 1}$ restricted to $\partial E$) (see \cite[Proposition 3.1]{CF1});

\smallskip
\item[\rm (ii)] For $\frac{n}{n - 1} \le p < \infty$,
$|\ban{\FF\cdot \nu, \cdot}_{\partial E}|(B) = 0$
for any Borel set $B \subset \partial E$
with $\sigma$--finite $\Haus{n - p'}$ measure
(see \cite[Theorem 3.2(i)]{Silhavy1}).
\end{itemize}
This characterization can be used to find classes of vector fields for which the normal trace
can be represented by a measure.
An important observation is that, for a constant vector field $\FF \equiv \vv \in \R^{n}$,
$$
\ban{\vv \cdot \nu, \cdot}_{\partial E}
= - \div(\chi_{E} \vv) = - \sum_{j = 1}^{n} v_{j} D_{y_{j}} \chi_{E}.
$$
Thus, in order that $\sum_{j = 1}^{n} v_{j} D_{y_{j}} \chi_{E}$ is a measure,
it is not necessary to assume that $\chi_E$ is a $BV$ function,
since  cancellations could be possible so that the previous sum could still be a measure.
Indeed, such an example has been  constructed (see \cite[Remark 4.14]{cct}) for a set $E \subset \R^{2}$
without finite perimeter
and a vector field $\FF \in \mathcal{DM}^{p}(\R^{2})$ for any $p \in [1,\infty]$
such that the normal trace of $\FF$ is a measure on $\partial E$.

In general,  even for an open set with smooth boundary,
$\ban{\FF \cdot \nu, \cdot}_{\partial E} \neq \ban{\FF \cdot \nu, \cdot}_{\partial \overline{E}}$,
since  the Radon measure $\div \FF$ in \eqref{normal trace def} is {\it sensitive} to small sets and
is not absolutely continuous with respect to the Lebesgue  measure in general.

Finally, we remark that a Gauss-Green formula for $\DM^p$ fields and
extended divergence-measure fields ({\it i.e.}, $\FF$ is a vector-valued measure whose divergence
is a Radon measure) was first obtained in Chen-Frid \cite{CF2} for Lipschitz domains.
In \v{S}ilhav\'y \cite{Silhavy2}, a Gauss-Green formula for extended divergence-measure
fields was shown to be also held over general open sets.
A formula for the normal trace distribution is given in
\cite[Theorem 3.1, (3.2)]{CF2} and \cite[Theorem 2.4, (2.5)]{Silhavy2} as the limit of
averages over the neighborhoods of the boundary.
In \cite{cct}, the normal trace is presented as the limit of classical
normal traces over smooth approximations of the domain. Roughly speaking,
the approach in \cite{cct} is to differentiate under the integral sign in the formulas
\cite[Theorem 3.1, (3.2)]{CF2} and \cite[Theorem 2.4, (2.5)]{Silhavy2}
in order to represent the normal trace as the limit of boundary
integrals ({\it i.e.}, integrals of the classical
normal traces $\FF\cdot\nu$ over appropriate smooth approximations of the domain).

\medskip
\section*{Entropy Solutions, Hyperbolic Conservation Laws, and $\DM$ Fields}

One of the main issues in the theory of hyperbolic conservation laws \eqref{codo}
is to study the behavior
of entropy solutions determined by the Lax entropy inequality \eqref{6-a}
to explore to the fullest extent
possible questions relating to large-time behavior, uniqueness,
stability, structure, and traces of entropy solutions, with neither
specific reference to any particular method for constructing the
solutions nor additional regularity assumptions.

It is clear that understanding more properties of
$\DM$ fields can advance our understanding of the
behavior of entropy solutions for hyperbolic conservation laws and
other related nonlinear equations by selecting appropriate entropy
pairs.
Successful examples
include
the stability of Riemann solutions,
which may contain rarefaction waves, contact discontinuities, and/or
vacuum states, in the class of entropy solutions of the Euler
equations for gas dynamics;
the decay of
periodic entropy solutions;
the initial and boundary layer problems;
the initial-boundary value problems;
and
the structure of entropy solutions of
nonlinear hyperbolic conservation laws.
See \cite{Chen2,CF1,CF2,ChenTorres,dafermos2010hyperbolic,Vas} and the references therein.

Further connections and applications of $\DM$ fields include
the solvability of the vector field $\FF$ for the equation: $\div \FF=\mu$ for given $\mu$,
image processing via the dual of $BV$,
and the analysis of minimal surfaces over weakly-regular
domains\footnote{See Meyer, Y.: {\it Oscillating Patterns in Image Processing and Nonlinear Evolution Equations},
AMS: Providence, RI, 2001;  Phuc, N. C., Torres, M.:
Characterizations of signed measures in the dual of $BV$ and related isometric isomorphisms,
{\it Ann. Sc. Norm. Super. Pisa Cl. Sci.} (5), 17(1): 385--417, 2017; Leonardi, G.~P., Saracco, G.:
Rigidity and trace properties of divergence-measure vector fields,
{\it Adv. Calc. Var.} 2021 (to appear), and the references therein.}.

Moreover, the $\DM$ theory
is useful for the developments of new techniques and tools for Entropy Analysis, Measure-Theoretic
Analysis, Partial Differential Equations, Free Boundary Problems,
Calculus of Variations, Geometric Analysis, and related areas, which involve the solutions
with discontinuities, singularities, among others.

\bigskip
\medskip
\noindent
{\bf Acknowledgements}.  The authors would like to Professors John M. Ball and John F. Toland,
as well as
the anonymous referees, for constructive suggestions to improve the presentation of this article.
The research of Gui-Qiang G. Chen was supported in part by
the UK Engineering and Physical Sciences Research Council Award
EP/L015811/1, and the Royal Society--Wolfson Research Merit Award (UK).
The research of Monica Torres was supported in part by the Simons Foundation Award No. 524190
and by the National Science Foundation Grant 1813695.

\smallskip
\medskip
Professor Gui-Qiang G. Chen is Statutory Professor in the Analysis of Partial Differential Equations and Director of the Oxford Centre for Nonlinear Partial
Differential Equations (OxPDE) at the Mathematical Institute of the University of Oxford, where he is also Professorial Fellow of Keble College.
$\,$ Email address: chengq@maths.ox.ac.uk

\bigskip
\smallskip
Professor Monica Torres is Professor of Mathematics at the Department of Mathematics of Purdue University.
$\,$ Email address: torresm@purdue.edu

\end{document}